\documentclass{amsart}
\usepackage{newlfont,amsfonts,amssymb,
amsmath,amsthm,amsgen,amscd,wasysym,comment}
\newcommand{\1}{\mathbb{I}}
\newcommand{\J}{\mathbb{J}}
\newcommand{\heta}{\hat{\theta}}
\newcommand{\ups}{\Upsilon}
\newcommand{\e}{\mathrm{e}}
\renewcommand{\d}{{\mathrm d}}
\newcommand{\bcase}{\begin{case}}
\newcommand{\ecase}{\end{case}}

\newcommand{\bclaim}{\begin{claim}}
\newcommand{\eclaim}{\end{claim}}

\newcommand{\bstep}{\begin{step}}
\newcommand{\estep}{\end{step}}

\newcommand{\bhlem}{\begin{hlem}}
\newcommand{\ehlem}{\end{hlem}}

\newcommand{\bleer}{\begin{leer}}
\newcommand{\eleer}{\end{leer}}
\newcommand{\bde}{\begin{de}}
\newcommand{\ede}{\end{de}}

\newcommand{\bs}{\begin{satz}}
\newcommand{\es}{\end{satz}}
\newcommand{\btheo}{\begin{theo}}
\newcommand{\etheo}{\end{theo}}
\newcommand{\bfolg}{\begin{folg}}
\newcommand{\efolg}{\end{folg}}
\newcommand{\blem}{\begin{lem}}
\newcommand{\elem}{\end{lem}}
\newcommand{\bnote}{\begin{note}}
\newcommand{\enote}{\end{note}}
\newcommand{\bprf}{\begin{proof}}
\newcommand{\eprf}{\end{proof}}
\newcommand{\bd}{\begin{displaymath}}
\newcommand{\ed}{\end{displaymath}}
\newcommand{\be}{\begin{eqnarray*}}
\newcommand{\ee}{\end{eqnarray*}}
\newcommand{\eeqa}{\end{eqnarray}}
\newcommand{\beqa}{\begin{eqnarray}}
\newcommand{\bi}{\begin{itemize}}
\newcommand{\ei}{\end{itemize}}
\newcommand{\bnum}{\begin{enumerate}}
\newcommand{\enum}{\end{enumerate}}
\newcommand{\la}{\langle}
\newcommand{\ra}{\rangle}

\newcommand{\ve}{\varepsilon}

\newcommand{\beq}{\begin{equation}}
\newcommand{\eeq}{\end{equation}}
\newcommand{\einhalb}{\frac{1}{2}}
\newcommand{\rr}{\mathbb{R}}

\newcommand{\vf}{\varphi}
\newcommand{\earr}{\end{array}\]}
\newcommand{\barr}{\[\begin{array}}
\newcommand{\bvec}{\left(\begin{array}{c}}
\newcommand{\evec}{\end{array}\right)}

\newcommand{\sumj}{\sum_{j=1}^n}

\newcommand{\lag}{\mathfrak{g}}

\newcommand{\lap}{\mathfrak{p}}

\newcommand{\+}{\oplus}

\newcommand{\laso}{\mathfrak{so}}

\newcommand{\w}{\omega}

\newcommand{\s}{\sigma}
\newcommand{\del}{\partial}
\newcommand{\ddt}{\frac{\partial}{\partial t}}

\newcommand{\ddu}{\frac{\partial}{\partial u}}

\newcommand{\bbem}{\begin{bem}}
\newcommand{\ebem}{\end{bem}}
\newcommand{\bbez}{\begin{bez}}
\newcommand{\ebez}{\end{bez}}
\newcommand{\bbsp}{\begin{bsp}}
\newcommand{\ebsp}{\end{bsp}}

\newcommand{\D}{\Delta}

\newcommand{\W}{\Omega}

\newcommand{\trace}{\mathsf{tr }}
\newcommand{\grad}{\mathsf{grad}}
\newcommand{\wt}{\widetilde}
\newcommand{\tnab}{\widetilde{\nabla}}
\newcommand{\tem}{\widetilde{M}}
\newcommand{\tg}{\widetilde{g}}

\newcommand{\ro}{\mathsf{P}}

\newcommand{\inter}{\makebox[7pt]{\rule{6pt}{.3pt}\rule{.3pt}{5pt}}\,}

\theoremstyle{definition}
\newtheorem{de}{Definition}
\newtheorem{bem}{Remark}
\newtheorem{bez}{Notation}
\newtheorem{bsp}{Example}
\theoremstyle{plain}
\newtheorem{lem}{Lemma}
\newtheorem{satz}{Proposition}
\newtheorem{folg}{Corollary}
\newtheorem{theo}{Theorem}

\begin{document}

\bibliographystyle{abbrv}


\title{Conformal structures with  $\mathrm{G}_{2(2)}$-ambient metrics}
\author{Thomas Leistner}\address[Leistner]{School of Mathematical Sciences, University of Adelaide, SA 5005, Australia} \email{thomas.leistner@adelaide.edu.au}
\author{Pawe\l~ Nurowski} \address[Nurowski]{Instytut Fizyki Teoretycznej,
Uniwersytet Warszawski\\ ul. Ho\.za 69, 00-681 Warszawa, Poland $\mathrm{and} $
Instytut Matematyczny PAN, ul. Sniadeckich 8, 
00-956 Warszawa, Poland}\email{nurowski@fuw.edu.pl}
%
\date{\today}
 \begin{abstract}
We present conformal structures in signature $(3,2)$ for which 
the holonomy of the Fefferman-Graham ambient metric is equal to the
non-compact exceptional Lie group $\mathrm{G}_{2(2)}$. We write down the
resulting 8-parameter family of $\mathrm{G}_{2(2)}$-metrics in dimension seven explicitly in an appropriately chosen coordinate
system on the ambient space.  
\\[.3cm]
\end{abstract}

\maketitle

\section{Introduction}
In Ref.~\cite{nurowski04} the second author constructed a conformal class
$[g_F]$ of $(3,2)$-signature metrics on every 5-manifold $M$ equipped with a
nonintegrable 2-distribution $S$ whose small growth vector is 
$(2,3,5)$. Such distributions in dimension five are called $(2,3,5)$-distributions. According to the classical results of Elie Cartan
\cite{cartan10} and David
Hilbert \cite{hilbert12}, $(2,3,5)$-distributions are in one to one
correspondence with ODE's of the form 
$$z'=F(x,y,y',y'',z), \quad\quad F_{y''y''}\neq 0,$$
for two real functions $y=y(x)$, $z=z(x)$ of one real variable
$x$. In particular, given $F$, the distribution is defined by 
$$S=\mathrm {Span}(\partial_{y''},\partial_x+y'\partial_y+y''\partial_{y'}+F\partial_z),$$ 
and the metrics $[g_F]$ are defined in terms of 
the function $F$ and its apropriate derivatives as in  \cite{nurowski04}. A construction of these metrics can also be found in \cite{cap-sagerschnig09}.

By the construction the metrics $[g_F]$ have reduced conformal
holonomy $H$. The group $H$ is contained in the noncompact form of the exceptional Lie 
group $\mathrm{G}_{2(2)}\subset \mathrm{SO}(4,3)$ \cite{nurowski04}. Moreover, it
was recently shown in \cite{hammerl-sagerschnig09}  that every conformal
class $[g]$ of $(3,2)$
signature metrics, whose conformal holonomy $H$ is contained in
$\mathrm{G}_{2(2)}\subset\mathrm{SO}(4,3)$ must be locally conformally equivenlent
to one of the structures $[g_F]$. 

Since metrics $[g_F]$ include all the conformal metrics with the
exceptional conformal holonomy $\mathrm{G}_{2(2)}$, it is interesting to ask
about the properties of their Fefferman-Graham ambient metrics 
$\tg_F$. 
The Fefferman-Graham ambient metric  of a conformal manifold $(M,[g])$ is a Ricci-flat metric on a neighbourhood in $\cal C\times \rr$  of the cone
\[\cal C=\{ ( g_p, p ) \mid p\in M, g\in [g]\},\]
which encodes the conformal structure (see \cite{fefferman/graham85, fefferman-graham07} or our Section \ref{ambientsec} for details).
In Ref.~\cite{nurowski07} properties of the ambient metric were studied for the conformal class 
 of metrics $[g_F]$ with $F=(y'')^2+a_0+a_1 y'+a_2 (y')^2+a_3 (y')^3+a_4
(y')^4+a_5 (y')^5 +a_6 (y')^6+b z$, where $a_\mu$, $\mu=0,\dots,6$, and $b$ are
real constants. The results about properties of $[g_F]$ with this $F$
from \cite{nurowski07} are strengthened in the present paper. 
Our aim here is to prove
\btheo\label{maintheo}
Let $(M,[g_F])$ be a conformal structure associated with a
$(2,3,5)$-distribution defined by a function \[F=(y'')^2+ a_0+a_1 y'+a_2 (y')^2+a_3 (y')^3+a_4
(y')^4+a_5 (y')^5 +a_6 (y')^6+b z,\] with $a_\mu$, $\mu=0,\dots,6$, and $b$ being 
real constants. This 8-parameter family of conformal structures has
the following properties:
\begin{enumerate}
\item \label{t3} For each value of the parameters $a_\mu$ and $b$ there exists a
  metric $g_F$ in the class $[g_F]$ and 
  ambient coordinates $(t,u)$, in which the  
Fefferman-Graham ambient metric $\tg_F$ for $g_F$ is given
explicitely by 
\[
\tg_F=-2\d t\d u+t^2 g_F-2 tu \ro-u^2 B.\]
Here $\ro$ and $B$ are the respective Schouten and Bach tensors for
$g_F$.
\item \label{t4} If at least one of $a_3$, $a_4$, $a_5$ or $a_6$ is not zero the
  metric $\tg$ has the full exceptional group $\mathrm{G}_{2(2)}$ as its 
pseudo-Riemannian holonomy. 
\end{enumerate}
\etheo
In the following proof of this theorem we will refer to statements that are proven in the paper.
\bprf
When proving the theorem we fix the metric $g_F\in [g_F]$ given as in formula \eqref{gF} on page \pageref{gF}.
In order to prove \eqref{t3} of the theorem, recall that 
the first term in the Fefferman-Graham ambient metric expansion is $(\mu_1)_{ab}=2\ro_{ab}$, where $\ro_{ab}$ denotes the Schouten tensor of $g_F$ (see \cite{fefferman-graham07} or our Section \ref{ambientsec} for details). 
The second term in the ambient metric expansion is given by
$(\mu_2)_{ab}=-B_{ab}+\ro_{ak}\ro^{k}_{\ b}$, where $B_{ab}$ is the Bach tensor of $g_F$. 
Using  the formulae for the Schouten tensor $\ro$ and the Bach tensor $B$ of $g_F$ given in the appendix, it follows that $\ro_{ak}\ro^{k}_{\ b}=0$. Again, using the formulae in the appendix,  one checks that
the metric
\[-2\d u\d t + t^2 g_F - 2tu \ro - u^2 B\]
is Ricci flat. Since $g_F$ is real analytic, the uniqueness of the ambient metric in the analytic category in odd dimensions \cite{fefferman/graham85, fefferman-graham07} implies that this is the ambient metric $\tg_F$ for $[g_F]$. This proves \eqref{t3} of the  theorem.

Now we prove \eqref{t4} of the theorem. 
Using the formulae for  $g_F$ given in the appendix, one verifies that
$\tg_F$ admits a parallel spinor which is not null (see 
Proposition  \ref{solution}).  Hence, the holonomy $H$ of $\tg_F$ is contained in $\mathrm{G}_{2(2)}$.
Now we have to verify that $H$ is equal to $\mathrm{G}_{2(2)}$.

First note that the naive approach of calculating the curvature of
$\tg_F$ and of showing that it generates $\mathrm{G}_{2(2)}$ does not
work since the curvature of $\tg_F$ is highly degenerated. In the best
case, the Riemann tensor of $\tg_F$ mapping $\Lambda^2\rr^7$ to
$\mathrm{G}_{2(2)} $ has rank four. On the other hand, obtaining the
full set of the first (or higher) derivatives of the curvature for the
general $F$ from the theorem is beyond our calculational skills. Thus we
have to use more subtle arguments. They are as follows: 

In order to verify that $H$ is equal to $\mathrm{G}_{2(2)}$, first assume that $H$ acts irreducibly on $\rr^{4,3}$. By Berger's list of irreducible holonomy groups of non-symmetric pseudo-Riemannian manifolds \cite{berger55}, which contains only $\mathrm{G}_{2(2)}$ in dimension $7$ (see also \cite{baumkath99} for the corresponding list of groups admitting invariant spinors), $\tg_F$ must be locally symmetric if $H\not= \mathrm{G}_{2(2)}$. 
This can be excluded by a direct calculation of derivatives of the 
curvature. For example, we verified by a direct calculation that 
$\tnab_1 \widetilde{R}_{1212} \neq 0$, where the indices refer
 to the orthonormal coframe $\xi^0, \ldots , \xi^6$ given on page \pageref{onframe}. 

Hence, if $H$ is not equal to $\mathrm{G}_{2(2)}$, it must admit an invariant subspace $V\subset \rr^{4,3}$. 
The exclusion of this situation will be based on the following two
pairs of statements. The first describes the relation between the
geometry of the ambient metric and the existence of certain metrics in
a conformal class. Let $(\tem, \tg)$ be the ambient metric for a
conformal class $[g]$. Then the following holds:
\begin{enumerate}
\item[(A)] If $(\tem, \tg)$ admits a parallel line bundle, then, on an open dense set of $M$, every metric in $[g]$ is locally conformal to an Einstein metric (see Theorem \ref{einsteintheo}). 
\item[(B)] If $(\tem, \tg)$ admits a parallel  bundle of totally null $2$-planes, 
then, on an open dense set of $M$, every metric in  $ [g]$ is 
locally conformal to a metric 
 $g$  which admits a parallel null line bundle $L$ such that $L^\bot\inter Ric^g=0$ (see Theorem \ref{theo-null}).
\end{enumerate}
The second pair of statements excludes the existence of certain metrics in the conformal class $[g_F]$ under assumptions on $F$.  In Theorems \ref{theorem0} and \ref{theorem01} we prove: If at least one of the constants $a_3$, $a_4$, $a_5$ or $a_6$ is
  not equal to zero, then
\begin{enumerate}
\item[(C)]\label{t1} the class $[g_F]$ does not contain a local
  Einstein metric (see Theorem \ref{theorem0}). 
\item[(D)] \label{t2a} the class $[g_F]$ does not contain a local
  metric $g$ that admits a $\nabla^g$-parallel null line $L$ whose
  Ricci tensor is annihilated by $L^\bot$ (see Theorem \ref{theorem01}).
\end{enumerate}

Now we assume that the   subspace $V$, which is invariant under the holonomy $H$ of $\tg_F$, is non-degenerate, i.e. $V\cap V^\bot=\{0\}$. By the local de Rham decomposition theorem, this implies that the ambient metric splits locally as a pseudo-Riemannian product metric, $\tg_F=g_1+g_2$. Since $\tg_F$ is Ricci flat, both, $g_1$ and $g_2$ have to be Ricci-flat. Since one of them is a metric in dimension $\le 3$, its Ricci-flatness implies that it is flat. In this case $\tg_F$ would admit at least one parallel vector field.
Now we use statement (A) that the existence of a parallel vector field for the ambient metric implies that, on an dense open set, $g_F$ is locally  conformal 
 to an Einstein metric. Under the assumptions on $F$, statement (C)
 gives the contradiction. Hence, $\tg_F$ does not admit a non-degenerate invariant subspace under the holonomy representation. 

Now assume that the $H$-invariant vector space $V$ is degenerate,
i.e. $W:=V\cap V^\bot$ is a non-trivial totally null space. Then the dimension of $V$ has to be  $\le 3$. 
The first case is that $W$ is one dimensional, i.e. that the ambient metric admits a parallel null line bundle.
Again by (A), $g_F$ must be locally  conformally
Ricci flat on an dense open set, which is again in 
contradiction to the statement (C).

Now, if $W$ is a null $2$-plane, statement (B) shows that locally there is a metric $g\in[g_F]$ with a $\nabla^g$-parallel null line $L$ and with $Ric^g(Y,.)=0$ for all $Y\in L^\bot$. By the assumptions on $F$ this contradicts the statement (D) that for the given $F$'s this is not possible.

Finally, we assume that $W$ is maximally null, i.e. three-dimensional. In this case there exists a pure null spinor which scales  under $H$ \cite{kath98g2,kopczynski97} and thus defines a line of spinors that is parallel for $\tg_F$. This means that we are in the situation where we have the parallel non null spinor defining $\mathrm{G}_{2(2)}$ and a parallel null line of spinors. 
Now we use the fact (which is proven in Lemma \ref{spinlemma}) that in this situation there exists a parallel line of vectors for $\tg_F$. Again, by 
(A)
this is in contradiction with (C) that $g_F$ is not conformally Einstein. This completes the proof.
\eprf

We emphasize that as a byproduct of this theorem we get explicit
formulae, in a coordinate system $(t,u,x,y,y',y'',z)$, for an
8-parameter family of strictly $\mathrm{G}_{2(2)}$ metrics in dimension
seven. These metrics have signature $(4,3)$ and are explicitely given
in formula \eqref{theoambient}. In Proposition \ref{solution} we also give the explicit expressions for
a parallel spinor $\psi$ for these metrics and furthermore the
explicit expressions for the corresponding closed and coclosed 
threeform $\w$ defining
the $\mathrm{G}_{2(2)}$ structure on the ambient space $(\widetilde{M},\widetilde{g})$.

%

\subsection*{Acknowledgements} We wish to thank Matthias Hammerl for
reading the preliminary version of this article and pointing out to us
a gap in our reasoning, which we eliminated in the final version. We also acknowledge helpful 
discussions with Andrzej Trautman, who provided us with a very useful 
representation of the real Clifford agebra Cl(4,3), which we use in Section \ref{spinsec}.  
This  work was supported by the SFB 676 of the German Research
Foundation and by the first author's Start-Up-Grant of the Faculty of
Engineering, Computer and Mathematical Sciences of the University of
Adelaide. The second author wishes to thank Rod Gover, Jerzy Lewandowski, and
the Mathematical Institute of Polish Academy of Sciences for making
possible his trip to the Southern Hemisphere, where the final version of this paper
was prepared.

\section{The ambient metric of an odd-dimensional  conformal structure}\label{ambientsec}
An important  tool in conformal geometry is the  {\em Fefferman-Graham ambient metric} (see \cite{fefferman/graham85} and \cite{fefferman-graham07} for the following). 
For a conformal class $[g]$ in signature $(p,q)$ on an
$n=(p+q)$-dimensional manifold $M$ one considers the cone 
\[\cal C=\{ ( g_p, p ) \mid p\in M, g\in [g]\}.\] We denote by $\pi:\cal  C\to M$ the canonical projection and by $\pi_*:T\cal C\to TM$ its differential. $\cal C$ is equipped with an obvious $\rr_+$-action 
$\vf_t(g_p,p)=(t^2g_p,p)$ and with the 
{\em tautological tensor} $G$ defined by
\[G_{(g_p,p)} (U,V):= g_p \left(\pi_*(U), \pi_*(U)\right).\]
The $\rr_+$-action extends to $\cal C\times \rr$. Now, the 
 the \emph{ambient space} $\tem$ with {\em ambient metric} $\tg$ is defined by the following properties:
 \bnum
 \item
 $\tem$ is an invariant neighbourhood of $\cal C$  in $\cal C\times
   \rr$ under the $\rr_+$ action.
 \item 
  $\tg$ is a smooth metric of signature $(p+1,q+1)$ on $\tem$ that is
   homogeneous of degree two with respect to the $\rr_+$-action and such that its pullback by $\iota:\cal C\to \tem$ gives the tautological tensor $G$, i.e. $\iota^*\tg=G$.
  \item The Ricci tensor $\widetilde{Ric}$ of $(\tem,\tg)$ is zero.
 \enum
 In the following we are only interested in the case where $M$ is odd-dimensional. In this case, Fefferman and Graham proved the following result.
 \btheo[\cite{fefferman/graham85,fefferman-graham07} and \cite{kichenassamy04}]
Let $(M,[g])$ be a real analytic manifold $M$ of {\em odd} dimension $n> 2$ equipped with a  conformal structure defined by a real analytic semi-Riemannian metric $g$. 
Then there exists an ambient space $(\wt{M},\wt{g})$ with real analytic Ricci-flat metric $\wt{g}$.
The ambient space 
 is unique modulo diffeomorphisms that restrict to the 
identity along ${\mathcal Q}\subset \wt{M}$ and commute with the $\rr_+$-action.
 \label{fg}
\etheo
By the uniqueness of the ambient metric, its pseudo-Riemannian
holonomy is an invariant object of the conformal class: Two different
ambient metrics corresponding to two different metrics from the same
conformal class are related by a diffeomorphism, $\tg_1=\psi^*\tg_2$,
and hence, are  isometric. Since the holonomy group is a
pseudo-Riemannian invariant, the holonomy group of the ambient metric 
in odd dimensions is a conformal invariant. 

Every  metric $g$ in the conformal class $[g]$ defines an embedding
\begin{eqnarray}\label{iota}
\iota_g:M\ni p\mapsto (g_p,p)\in \cal C
\end{eqnarray}
and thus an identification of $\tem$ with $\rr_+\times M\times \rr$ via
\[ (t,p,\rho)\mapsto (t^2 g_p,\rho)\]
Using this identification and 
starting with a formal power series
\begin{equation} 
\tg = 2 \left( t \d\rho + \rho \d t \right) \d t + t^2
\left( g + \sum_{k=1}^\infty \rho^k\mu_k \right)\label{afg}
\end{equation}
with 
 certain symmetric $(2,0)$ tensors  $\mu_k$ on $M$,
 Fefferman and Graham showed that if $n$ is \emph{odd},
the Ricci-flatness of the ambient metric gives equations for 
$ \mu_1,\mu_2, \ldots $ that can be solved. 
However, the $\mu_k$ have been determined only for small $k$  or for all $k$ but very special conformal classes. 
For example, in general one finds that $\mu_1=2\ro^g$, where $\ro$ denotes the Schouten tensor of $g$, and that
\begin{equation}\label{mu2}\mu_2= -B^g + \trace ( \ro^g \otimes \ro^g)\end{equation}
with $B^g$ being the Bach tensor of $g$. Furthermore,
for an Einstein metric with $\ro^g= \Lambda g$ we have that $\mu_2=\Lambda^2 g$ and all other $\mu_i=0$, i.e. the power series in the ambient metric $\tg_E$ \emph{truncates}
at $k=2$.  
Further calculations of the ambient metric
 have been carried out for conformal 
classes that are related to Einstein spaces 
\cite{gover-leitner06}.   
However, if the metric $g$ is 
\emph{not conformally Einstein}, then, except for a few examples 
\cite{gover-leitner06,nurowski07,leistner-nurowski08}, no explicit formulae for 
$\mu_k$, $k>3$ are known.

For further convenience we change  the coordinate $\rho$ on $\tem$
to
 $u:=-\rho t$, i.e. $\d u=-t\d \rho - \rho \d t$. Then the ambient metric takes the form
\be
\tg&=& -2\d u\d t +t^2g -2 ut \ro^g +u^2 
\left(\mu_2 - \frac{u}{t} \mu_3 +\left(\frac{u}{t}\right)^2\mu_4- ...\ \right).
\ee

In particular, for a Ricci-flat metric, the ambient metric is given as  a special  Brinkmann wave,
\begin{equation}
\label{ricflat}
\tg= -2\d u \d t +t^2g,\end{equation}
admitting a parallel null vector field,
whereas for an Einstein metric with  $\ro=\Lambda g$ the ambient metric becomes 
\[
\tg= -2\d u\d t +\left( t^2 -2\Lambda ut +\Lambda^2u^2\right) g.
\]
This metric 
splits into a line and a cone. This becomes evident  in new coordinates $r=t-c u$ and $s=t+\Lambda u$  in which we have 
\begin{equation}
\label{cone}
\tg= \frac{1}{2\Lambda }(\d r^2 -\d s^2)+r^2 g.
\end{equation}

Now, let $\tg$ be the ambient metric for an arbitrary conformal class $[g]$ on an odd-dimensional manifold. We calculate the Levi-Civita connection $\tnab$ of $\tg$ 
along $\cal Q=\{u=0\}$ and obtain that the only non-vanishing terms are 
\begin{equation}\label{nabamb}
\left.\begin{array}{rcl}
\tnab_X\ddu&=& -\frac{1}{t}\ro^g(X)^*\\
\tnab_X Y&=&\nabla_X Y+ t\left( g(X,Y)\ddu -  \ro^g(X,Y)\ddt\right)\\
\tnab_X\ddt&=&\frac{1}{t} X\\
\end{array}\right\},\end{equation}
for $X,Y\in \Gamma (TM)$ and $\nabla$ being the Levi-Civita connection of $g$.

In the following we need the transformation of the Schouten tensor
under a conformal rescaling. Recall that if $\hat g=\e^{2\ups}g$ with
$\ups\in C^\infty (M)$ is a conformally changed metric, then the Schouten tensor $\hat \ro$ of $\hat g$ satisfies
\begin{equation}\label{schoutentrafo}
\hat \ro  =\ro  - \mathrm{Hess}^g(\ups) +\d \ups^2 -\einhalb \|\grad^g(\ups)\|^2_g g,
\end{equation}
where $\mathrm{Hess}(\ups)=g(\nabla \grad^g(\ups),.)$ denotes the Hessian of $\ups$. For brevity we will also write this relation as
\[
\hat\ro_{ab}=\ro_{ab}-\nabla_a\ups_b+\ups_a\ups_b -\einhalb \ups_c\ups^c g_{ab}.\]
Hence, $\hat g$ is an Einstein metric if and only if this quantity is a multiple of the metric $g$.
More explicitly, it holds that $\hat\ro_{ab}=\Lambda \hat g_{ab}$, where $\Lambda $ is a constant, if and only if
\[
\ro_{ab}-\nabla_a\ups_b+\ups_a\ups_b = \left( \einhalb \ups_c\ups^c +\Lambda \e^{2\ups}\right) g_{ab}
\]
By substituting $\ups =-\log(\s)$ for a non vanishing  function $\s$, we obtain that the metric $\hat g= \s^{-2}g=\e^{2\ups}g$ is Einstein if and only if there is an non vanishing function $\s$ such that
the symmetric tensor $\mathrm{Hess}(\s)+\s\ro$ is a multiple of $g$. Explicitly, we have
\begin{equation}
\label{alpha}
\nabla_a\s_b+\s\ro_{ab} = \s^{-1}\left( \Lambda +\einhalb \s_c\s^c\right)g_{ab}.
\end{equation}
Now we give a criterion for a metric being locally conformally Einstein in terms of the ambient metric.
\btheo\label{einsteintheo}
Let $M$ be an odd dimensional manifold equippped with a real analytic conformal class $[g]$. If the ambient space $(\tem,\tg)$ admits a line bundle $\cal L$ that is parallel with respect to the Levi-Civita connection of $\tg$, then on the connected components of an open dense subset $M_0$ in $M$, every metric in the conformal class $[g]$ is locally conformal to an Einstein metric $g_E$.

Furthermore, if the ambient metric on $\cal L$ is positive/negative/null, then the Einstein constant  of $g_E$ is negative/positive/null.
\etheo
\bprf
We start the proof with a lemma.
\blem
There is no open set  $U$ in $ \cal C$ such that $\cal L|_U\subset T\cal C|_U$.
\elem
\bprf
Assume that we have an open set $U\subset \cal C$ such that $\cal L|_U\subset T\cal C|_U$. Since $\cal L$ is parallel, by making $U$ smaller, such that it becomes simply connected, we can assume that there is a section 
$L\in \Gamma (\cal L|_U)$. We fix a metric $g\in[g]$ to obtain $(t,x)$-coordinates on $U$ and write
\[L=\alpha \del_t+ K,\]
with $K$ tangential to $M$.  Since $\cal L$ is parallel,  formula \eqref{nabamb} implies
\[
0\equiv \tg(\tnab_XL,\del_t)
=-
 \tg(L,\tnab_X\del_t)
 =
 \tfrac{1}{t}\tg(X,K) = tg(X,K),\]
for all $X\in TM$. This implies $K\equiv 0$ on $U$. Hence $L= \alpha\del_t$, but this contradicts $\tnab_X\del_t=\tfrac{1}{t}X$ for all $X\in TM$.
\eprf

This lemma implies that there is an open dense set $\cal C_0$ in $\cal C$ such that $\cal L|_{\cal C_0}\not\subset T\cal C|_{\cal C_0}$. For every point in $M_0:=\pi(\cal C_0)$ we have to verify the existence of a neighbourhood on which a metric in $[g]$ can be rescaled to an Einstein metric. The following lemma will be useful.
\blem
On every simply connected  open subset $U$ of the open and dense subset $\cal C_0$ in $\cal C$ there is a section $L\in \Gamma (\cal L|_U )$ such that $\tnab_Y L|_{U}=0$ for all $V\in T\cal C$.
\elem
\bprf
For every simply connected open set $U$ in $\cal C_0$ we find a section of $\cal L$ which, by fixing $g\in [g]$ and by the previous lemma, is of the form
\[L=a\del_t+K+\del_u\] with $K$ tangential to $M$.
 Since $\cal L$ is parallel, there is a 1-form $\Theta$ over $U$ such that $\tnab L=\Theta\otimes L$. We will show  that $\Theta $ is closed, which implies that  $L$ can be rescaled to a parallel vector field.
 The following calculations are over $U$. We get
\[
\Theta(\del_t) = -\tg(\tnab_{\del_t} L,\del_t)\equiv 0, \]
which implies that
\[K=\tfrac{1}{t}K_0\]
 for a $K_0\in \Gamma(TM)$. Indeed,  for every $X\in TM$ it is
\[
\tg(\tnab_{\del_t}L,X) = \tg(\tnab_{\del_t}K,X) =\tg ( [\del_t,K], X) +\tfrac{1}{t}\tg(K,X).\]
Now, as 
$\tg(\tnab_{\del_t}L,X)=\Theta(\del_t) \tg(K,X)=0$, this implies that $K$ satsifies the equation
\[ \left[\del_t,K\right] = -\tfrac{1}{t}K,\]
which yields $K=\tfrac{1}{t}K_0$ with $K_0\in \Gamma(TM)$.
Furthermore it is
\[
\Theta (X)= 
-\tg(\tnab_X L,\del_t)= \tg( L,\tnab_X\del_t)= \tfrac{1}{t}\tg(K,X)=g(K_0,X)\]
for $X\in TM$.
Hence, in order to show that $\Theta$ is closed we only have to check $\d\Theta (X,Y)=0$ for $X,Y\in TM$.
On the one hand we get that
\[
\tg(\tnab_X L, Y) = t \left( a g(X,Y) +t g(\nabla_XK_0, Y) -P^g(X,Y)\right),\]
and on the other that
\[
\tg(\tnab_X L, Y)= t \Theta(X) g(L,Y) =t \Theta^2(X,Y),\]
which shows that 
$g(\nabla_XK_0, Y)$ is symmetric in $X,Y\in TM$. But this implies that $\Theta=g(K_0,.) $ is closed.
Hence, on simply connected open sets $U\subset \cal C_0$ we get that $\Theta =\d f$ which implies that $\e^{-f}\cdot L $ 
is a parallel vector field on $U$.
\eprf
Now we conclude the proof of the theorem by fixing a metric $g$ in
$[g]$ and showing that it can be rescaled to an Einstein metric $g_E$ on simply connected open sets in $M$. By the lemmas, on simply connected open sets  in $\cal C_0$ we get a parallel vector field 
\[ L=\alpha \del_t+K+\s\del_u\in \Gamma(\cal L|_U),\]
with $\s\not=0$ and $K$ tangential to $M$.
Again, $ \tnab_{\del_t}L=0$ implies that $K=\tfrac{1}{t}K_0$ with $K_0\in \Gamma (TM)$, but also that $\d \s(\del_t)=\d \alpha (\del_t)\equiv 0$. For $X\in TM$, the equation $\tnab_XL=0$ implies
\begin{eqnarray}
\label{sigma}
\d \s (X)+g(X,K_0)&=&0\\
\label{a}
\d \alpha  (X)-P(X,K_0)&=&0\\
\label{grad}
g(\nabla_XK_0, Y)- \s P^g(X,Y) + \alpha g(X,Y)&=&0
\end{eqnarray}
The first equation shows that $K_0=-\grad^g(\s)$. Then the last equation shows 
\[
\mathrm{Hess} (\s)+\s P^g = \alpha g.\] 
But  this is equivalent to $\s^{-2} g$ being a local Einstein metric.
Note that \eqref{alpha} implies that 
\[\s\alpha = 
\left( \hat\Lambda +\einhalb g(K_0,K_0)\right)
\]
with the Einstein constant $\hat\Lambda$ of  $\s^{-2}g$.
But this implies that 
\[\tg(L,L) = -2\alpha \s +g(K_0,K_0)= -2\hat\Lambda, \]
which shows the relation between the Einstein constant and the line bundle being null, positive, or negative.
\eprf

\bbem
Using the formulae in \eqref{nabamb}, one can show that the connected component of the normal conformal Cartan connection for the conformal structure $[g]$ is contained in the holonomy of the ambient metric $\tg$. 
Furthermore,
when  the conformal class contains an Einstein metric, the truncation of the ambient metric in this case yields the equality of both holonomy groups (see  \cite{leistner05a} and \cite{leitner05}). 
This can be used to prove an analogue of  Theorem \ref{einsteintheo}
in terms of the normal conformal Cartan connection. This analogue
holds in any dimension and gives an equivalence between the existence
of a parallel line ${\mathcal L}$ in $(\tem,\tg)$ and an Einstein metric
$g_E$ in $[g]$.
As we will use here only one direction and only in odd dimensions, for the purpose of being self contained, we did prove the Theorem without referring to the normal conformal Cartan connection and without using tractor calculus. For further results relating the ambient metric and tractor calculus, see \cite{cap-gover03, leistner05a, armstrong-leistner07}.
\ebem
 
Now we will describe the case where the Levi-Civita connection of the
ambient metric admits an invariant null $2$-plane\footnote{The case
  where the $2$-plane is non-degenerate implies that the conformal class
  contains a product of Einstein metrics with related Einstein
  constants (see \cite{armstrong07conf} for Riemannian conformal
  classes and the unpublished parts in \cite{leitner04} for arbitrary
  signature).}. We will deal with a bit more general situation than
needed for our pourposes, i.e. with the case when the $2$-plane is totally
null in arbitrary signature. The following theorem is a generalisation
to arbitrary signature of the corresponding  result from the Lorentzian domain, which was proved in \cite{leistner05a}.

\btheo\label{theo-null}
Let $(M,[g])$ be a pseudo-Riemannian real analytic conformal manifold of odd dimension $n>2$. If the holonomy group of the ambient metric admits a parallel null $2$-plane, then every  metric  in the conformal class  $[g]$ is locally conformally equivalent to a metric $g$ with the following two properties:
\begin{eqnarray}
\label{null-line}
&&\text{There is a null line $L\subset TM$ that is parallel for $\nabla^g$, and }
\\
&&Ric^g(Y)=0 \text{ for all }Y\in L^\bot.\label{null-ric}
\end{eqnarray}
\etheo
\bbem
Note that property \eqref{null-ric} is equivalent to the property that the image of $Ric^g:TM\to TM$ is contained in $L$.
In Lorentzian signature, this is equivalent to the image of $Ric^g$ being totally null. In higher signature, it is stronger.

Note also that such metrics have vanishing scalar curvature, and thus
$Ric$ is a constant multiple of the Schouten tensor. It also holds that $Ric$ satsifies property \eqref{null-ric} if and only if the Schouten satisfies  property \eqref{null-ric}.
\ebem
\bprf
The proof is based on the following Lemma.
\blem\label{lem1}
Let $(\tem=\cal C\times \rr , \tg)$ be the ambient space and let $\cal H$ be a bundle of parallel null $2$-planes on $(\tem, \tg)$.
Assume that there is a bundle of null lines $L$ over $M$ and a metric $g\in [g]$ defining the embedding $\iota_g:M\to \cal C$ such that 
\[
\cal H|_{\iota_g(M)} =\ (\iota_g)_*(L)\+ \rr \partial_u.\]
Then $L$ is parallel with 
respect to $\nabla^g$ and $Y\inter Ric^g =0$ for all $Y\in L^\bot$.
\elem
\bprf
Let $g\in [g]$ be the metric given in the assumptions and let $\tem = \rr_+\times \iota_g(M) \times \rr\ni (t,p,u)$.
Furthermore, let $V=a \del_u+K\in \Gamma (\cal H|_{\iota_g(M)})$ with $K\in \Gamma (L) $ a null vector. 
Then, along $\iota_g(M)\subset \cal C$, by formulae \eqref{nabamb},  we get
\[0=\tg(\tnab_X V, \del_u)=
-\tg(\tnab_X  \del_u,K)= tP^g(X,K),\]
for all $X\in TM$.
This shows that $L\inter P^g=0$. In particular, the image of $Ric^g$ and hence the image of $\ro^g$ lies in $L$. Furthermore,
\[\tnab_X V= X(a)\del_u +\frac{a}{t}\ro^g(X) + \nabla^g_XK+ tg(X,K)\del_u,\]
for all $X\in TM$.
Since $\cal H|_{\iota_g(M)}=L\+\rr\del_u$ is parallel, and since $\ro^g(X)\in L$, this implies that $L$ is parallel with respect to $\nabla^g$.\eprf
We will now show that the existence of a parallel totally null $2$-plane distribution on $(\tem,\tg)$ implies the existence of a metric in the conformal class and a null line bundle on $M$ satisfying the assumptions of Lemma  \ref{lem1}.

Let $\cal H$ be a totally  null $2$-plane bundle  that is parallel for the ambient metric. With $\cal H$ also $\cal H^\bot$ is parallel, but of rank $n$.
This implies
\begin{equation}\label{notsub}
\cal H|_{\cal C}\not\subset T\cal C\ \text{ and }\ \cal H^\bot|_{\cal C}\not\subset T\cal C.\end{equation}
In order to prove this, fix a metric $g\in [g]$. Since $\cal H$ has rank $2$, $\cal H|_{\cal C} \subset T\cal C$ would imply that there is a section $K\in \Gamma (\cal H|_{\cal C})$ that is tangential to $\iota_g(M)$. 
$\cal H$ being parallel then gives
\[ 0=\tg(\tnab_XK,\del_t) =- \tg(\tnab_X\del_t,K) = -tg(X,K),\]
for all $X$ in $TM$, which contradicts the non-degeneracy of the metric.
Then, since $\cal H$ is totally null, and thus $\cal H\subset \cal H^\bot$, property \eqref{notsub} follows.
For reasons of dimensions, this implies that
\begin{equation}
\label{line}
\cal L:=\cal H|_{\cal C}\cap T\cal C\text{ is a bundle of null lines over $\cal C$.}
\end{equation}
We will now prove some properties of $\cal L$ that will lead to the proof of the theorem.
\blem
Let $\cal X:=\cal L^\bot\cap T\cal C$ and $L\in \Gamma(\cal L)$. Then
\begin{equation}
\label{wichtig}
\tg(\tnab_U L,V)=0 \ \ \ \text{ for all }U,V\in \cal X.
\end{equation}
\elem
\bprf
We show that for $U\in \cal X$ the vector field  $\tnab_UL$ is not
only contained in $\cal H$, by $\cal H$ being parallel, but also in
the space tangential to the cone, and hence in $\cal L$. To this end we fix a metric $g\in [g]$ yielding $\cal C=\rr_+\times M$ and $T\cal C=\del_t^\bot$. Hence, 
for $L=a\del_t+K\in \cal L$ and $U=b\del_t+X\in \cal X$, with $K$ and
$X$ tangential to $M$ and orthogonal to each other, 
\[
\tg(\tnab_{U} L,\del_t)=
-\tg(L,\tnab_U\del_t) =-\tg(L,\tnab_X\del_t)= 
- \tfrac{1}{t}\tg(L,X)
=
- \tfrac{1}{t}\tg(X,K)
=0,
\]
since $g(K,X)=0$ and $\del_t^\bot= T\cal C$.
This shows that $\tnab_{U} L\in \cal H\cap T\cal C$ implying the relation \eqref{wichtig}.
\eprf

\blem
$\cal X=\cal L^\bot\cap T \cal C$ is an integrable distribution on $\cal C$.
\elem
\bprf
We fix a metric $g\in [g]$ and obtain $M\hookrightarrow \cal C$.
First note that $\cal H^\bot\subset \cal L^\bot$ and that $\del_t\in \cal L^\bot$. Furthermore $\cal H^\bot \cap \del_t=\{0\}$. Indeed, since $\cal H\not\subset T\cal C$, there is an element in $\cal H$ of the form $
\del_u+X+a\del_t$, which implies that $\del_t$ is not orthogonal to $\cal H$. Furthermore, since $\cal H^\bot\subset \cal L^\bot$, by relation \eqref{notsub}, the dimension of $\cal L^\bot\cap \cal T\cal C$ is $n$. 
Now, let $\del_t,Y_1, \ldots , Y_{n-1}$ with $Y_i\in \Gamma (\cal H^\bot\cap T\cal C)$ be basis for $\cal L^\bot\cap T\cal C$ of mutually orthogonal vector fields. Then, since $\cal H^\bot $ is parallel we have $\tnab_{\del_t}Y_i\in  \Gamma(\cal H^\bot)$. Hence,
\[\left[ \del_t,Y_i\right] = \tnab_{\del_t}Y_i -\tnab_{Y_i}\del_t = \tnab_{\del_t}Y_i -\tfrac{1}{t}Y_i\ \in\ \cal H^\bot\cap T \cal C.\]
On the other hand, again since $\cal H^\bot$ is parallel, we get $\tnab_{Y_i}Y_j\in \Gamma(\cal H^\bot)$ and thus $\left[Y_i,Y_j\right]\in \Gamma(\cal H^\bot)$. Furthermore
\[ \tg\left( \left[ Y_i,Y_j\right],\del_t\right)= \tg\left( \tnab_{Y_i}Y_j - \tnab_{Y_j}Y_i,\del_t\right)
=
- \tg\left( Y_j ,\tnab_{Y_i} \del_t\right) +  \tg\left( Y_i ,\tnab_{Y_j} \del_t\right)
=
0,
\]
because of \eqref{nabamb}.
This shows that also $\left[ Y_i,Y_j\right]\in \Gamma (\cal H\cap T\cal C)$. Hence, $\cal H^\bot\cap T\cal C$ but also
$\cal L^\bot\cap T\cal C$
are integrable.
\eprf

\blem 
Let $\pi:\cal C\to M$ and $\pi_*:T\cal \to TM$ be the canonical projection. Then 
\[L:=\pi_*(\cal L) \]
is a distribution of null lines on $M$ and $L^\bot=\pi_*(\cal X)$. Both distributions are integrable on $M$.
\elem
\bprf
We fix $g\in [g]$ to verify that $\pi_*(\cal L)\not=0$. Assume that $\cal L=\rr \del_t$. $\cal H$ being parallel then implies
that 
\[
\tnab_X\del_t=\frac{1}{t} X \in \cal H\ \ \text{ for all $X\in TM$.}\]
Since $n>2$, this contradicts to $\cal H$ being a $2$-plane bundle. Hence, $\pi_*(\cal L)$ is a null line bundle on $M$. Since $T\cal C=\del_t^\bot$ this implies that $L^\bot=\cal L^\bot\cap T\cal C=\cal X$. Since $\cal L$ and $\cal X$ are integrable on $\cal C$, $L$ and $L^\bot$ are integrable on $M$.
\eprf

\blem
For any null vector $K\in \Gamma ( L) $ we define the second fundamental form of $L^\bot$ by
\[\Pi^K(X,Y)
=g(\tnab_X K,Y)\ \ \text{ for }X,Y\in L^\bot.\]
Then $\Pi^K$ is symmetric and tensorial in  $K$. Furthermore, locally there is a metric in the conformal class such that
$\Pi^K(X,Y)$ has no trace.
\elem
\bprf
First we notice that $L\subset L^\bot$ implies that $\pi^K$ is tensorial in $K$: For $f K$ with a smooth function $f$ we get
\[
\Pi^{fk}(X,Y)= X(f) \tg(K,X)+ f \tg(\tnab_XK,Y)= f\Pi^K(X,Y).
\]
The integrability of $L^\bot$ implies that $\Pi^K$ is symmetric.
Now we define the trace of $\Pi^K$ as
\[H^K:=\sum_{i=1}^{n-2}\ve_i \Pi^K(E_i,E_i)\in C^\infty (M)\]
where $E_1,\ldots , E_{n-2}$ linearly independent in $L^\bot$ with $g(E_i,E_j)=\ve_i\delta_{ij}$. Since $K\inter \Pi^K=0$, this is independent of the chosen $E_i$'s.
Now we claim that there is a metric $\hat{g}=\e^{2\Upsilon}g\in [g]$ in the conformal class such that the corresponding function $\hat{H}^K$ is zero. To this end we notice that the transformation formula for $\hat{\Pi}^K$ is given by
\[
\hat{\Pi}^K(Y,V)= \hat{g}(\hat{\nabla}_YK,V)
=
\e^{2\Upsilon}\left(\Pi(Y,V)+\d \ups (K) g(V,Y)\right),
\]
for $Y,V\in L^\bot$. Hence,
\[\hat{H}^K=\e^{2\ups}\left( H^K+ (n-2)\d \ups(K) \right).\]
Now the differential equation 
\[\d\ups(K)=\frac{H^K}{n-2}\]
is an ODE along the flows of $K$ and as such
locally always has a solution. This ensures that we can chose $\hat{g}$ such that $\hat{H}^K\equiv 0$.
\eprf

 Finally, to conclude the proof, we fix this metric $g\in [g]$ for which $H^K\equiv 0$. Now let $L=a\del_t+K\in \Gamma(\cal L)$ be arbitrary. Then equation \eqref{wichtig}
 reads as
 \[0= \tg(\tnab_X L,Y)=at g(X,Y)+t^2 \Pi^K(X,Y),\]
 for all $X,Y\in L^\bot$.
 Taking the trace shows that $a\equiv 0$ on $M$. Hence, $\cal L
 =\iota_g(L)$, and thus $\cal H= \iota_g(L)\+\rr\del_u$. This means
 that the metric $g\in [g]$ and the null line bundle $L$ on $M$ satisfy the assumptions of Lemma \ref{lem1}. This concludes the proof of the theorem.
 \eprf





\section{$\mathrm{G}_{2(2)}$-conformal structures with truncated ambient metric}
As it was mentioned in the Introduction, in \cite{nurowski04} a  conformal structure $[g_F]$ in signature $(3,2)$ was introduced that originated from a first order ODE for two functions $y,z$ of one variable $x$. We will now describe this construction briefly. Every solution to the first order ODE
\be
z'=F(x,y,y',y'',z)\quad\quad\text{ 
with }\quad\quad F_{y''y''}\neq 0,\label{2m}\ee
 is a curve in the five-dimensional manifold $M$ parametrised by $(x,y,z,p=y',q=y'')$, on which the one-forms 
\begin{equation}
\w^1=\d z-F(x,y,p,q,z)\d x,\ \ 
\w^2=\d y-p\d x,\ \ 
\w^3=\d p-q\d x
\label{fhil}
\end{equation} vanish.
Two triples of such 1-forms on  $ \rr^5$, $(\w^1,\w^2,\w^3)$ and $(\hat{\w}^1,\hat{\w}^2,\hat{\w}^3)$, are considered to be equivalent, if there is a local diffeomorphism $\Phi$ of $\rr^5$ and a $\mathrm{GL}(3,\rr)$-valued function $A=(a_{ij})$ on the domain of $\Phi$ such that  $\Phi^*\hat{\w}^i=\sumj a_{ij}\w^j$.
 Cartan showed that an equivalence class of a triple of one-forms given by (\ref{2m}) with $F_{qq}\not=0$ corresponds to a Cartan connection $\w$ on a $14$-dimensional principle fibre bundle $\cal P$ over the five-manifold parametrised by $(x,y,z,p,q)$. This Cartan connection has values in the non-compact exceptional Lie algebra $\lag_{2(2)}$, and $\cal P$ is the bundle with structure group given by the $9$-dimensional parabolic $P:= \mathrm{G}_{2(2)}\cap B$, where $B$ is the isotropy group in $\mathrm{SO}(4,3)$ of a null line. The conformal structure on the 5-manifold is now constructed as follows:  Write the Cartan connection $\w$ as $\w=(\theta, \W)$, where $\W$ has values in the Lie algebra $\lap$ of $P$ and $\theta$ in the five-dimensional complement of $\lap$ in $\lag_{2(2)}$. Write $\theta=(\theta_1, \ldots , \theta_5)$ and $\Omega=(\W_1, \ldots , \W_9)$ and let $X_1, \ldots , X_5$ and $Y_1, \ldots , Y_9$ be the vector fields on $\cal P$ dual to $\theta_i$ and $\W_\mu
 $, respectively. The $Y_\mu$ are tangential to the fibres of $\cal P\rightarrow M$. Defining the bilinear form
 \[ G= 2\theta^1\theta^5 - 2\theta^2\theta^4+(\theta^3)^2\]
 on $\cal P$ we note that along the fibres $G$ is degenerate and merely  scales, i.e. 
 \[\cal L_{Y_\mu}G=  \lambda_\mu G\]
 for some functions $\lambda_\mu$.  Hence, $G$ projects to a conformal class of metrics  $[g_F]$ of signature $(+++--)$ on $M$.  This means  that the normal conformal Cartan connection for $[g_F]$ reduces (in the Cartan sense) to $\mathrm{G}_{2(2)}$. Hence, the conformal holonomy of $[g_F]$ is contained in this group. Of course, this inclusion might be proper.
 
 \bbem
 The conformal structure given by $F$ is an example of a {\em conformal Cartan reduction} (see for example \cite{alt08}). The normal
conformal Cartan connection  of $[g_F]$
 reduces to a Cartan connection with values in  the Lie algebra $\lag_{2(2)}\subset \mathfrak{so}(4,3)$. In this way it  defines a parabolic geometry of type $( P,\lag_{2(2)})$, where $P$ is the parabolic subgroup given by the stabiliser in $\mathrm{G}_{2(2)}$ of a null line. 
This situation is exceptional in the sense that a reduction of a Cartan connection to a semisimple subalgebra $\lag\subsetneq \laso(p+1,q+1)$ whose intersection with the stabiliser of a null line is parabolic imposes very strong algebraic restrictions on $\lag$ and the parabolic subalgebra, as recently shown in \cite{doubrov-slovak08}. For conformal geometry, only two cases arise: the one of $\mathfrak{g}_{2(2)}$, which, by the result in \cite{hammerl-sagerschnig09}, is given by the above  construction,  and the one of $\laso (4,3)\subset \laso (4,4)$ described in \cite{bryant06}. 
 \ebem
 
Then, in \cite{nurowski07}, the following remarkable feature of $[g_F]$ was noticed.
\bs
There exist functions  $F$ such that  the ambient metric of a $g_F\in
[g_F]$ truncates after terms of second order, i.e.
\begin{equation}\wt{g}_F= -2\d t\d u+t^2 g_F-2ut\ro+u^2\beta,\label{ambienttrunc}\end{equation}
with $\ro$ the Schouten tensor of $g_F$ and $\beta=\mu_2$ defined as in Eq. (\ref{mu2})
\es

Examples of such $F$'s given in \cite{nurowski07} include $F=F(q)$ and $F=
q^2+ \sum_{i=0}^6 a_i p^i+b z
$. 
The proof is based on the form and the uniqueness of the ambient metric in odd dimensions proved in \cite{fefferman/graham85,fefferman-graham07} and the observation, that the metric (\ref{ambienttrunc}) is Ricci-flat.

This concise form of the ambient metric makes it possible
 to study the relation between the conformal holonomy  and the holonomy of the ambient metric. 
This is done by distinguishing  two situations: the first, when the
conformal class contains an Einstein metric, and the second, when it does not contain an Einstein metric.
Also in \cite{nurowski07}
several examples of such conformal structures depending on the function $F$ in (\ref{2m}) with $F_{qq}\not=0$  were considered.  On the one hand it was shown that for  
$F=F(q)$
the conformal class given by $F$ contains a Ricci flat metric. We have seen that for a conformal class that contains a Ricci flat metric, 
the ambient metric is 
a special Brinkmann metric, $\tg=-2\d u \d t +t^2g$,  and that the holonomy of the ambient metric is the same as the holonomy of the conformal Cartan connection. Based on the result in \cite{nurowski07} we obtain:
\bs
Let $[g_F]$ be a conformal class where $F=F(q)$ with $F_{qq}\not=0$. Then $[g_F]$ contains a Ricci flat metric, the ambient metric for $[g_F]$ is 
\[\tg_F= -2\d u \d t +t^2g_F\]
as in equations (\ref{ricflat}), 
the holonomy of the ambient metric is equal to the conformal holonomy and contained in the  eight-dimensional stabiliser in $\mathrm{G}_{2(2)}$ of a null vector.
\es 
This shows that the ambient metric of conformal classes $g_{F(q)}$ are $\mathrm{G}_{2(2)}$-metrics that admit a parallel null vector field, and thus can be considered as $\mathrm{G}_{2(2)}$-Brinkmann waves.

Furthermore, in \cite{nurowski07} a conformal structure $[g_F]$ in signature $(3,2)$ was introduced that still has an ambient metric in the truncated form \eqref{ambienttrunc} but does {\em not} contain an Einstein metric. This is defined by 
\begin{equation}\label{eff}
F=q^2+ \sum_{i=0}^6 a_i p^i+b z.
\end{equation}
Explicitly, 
\begin{equation}\label{gF}g_F=2\theta^1\theta^5-2\theta^2\theta^4+(\theta^3)^2,\end{equation}
where the co-frames $\theta^i$ are given by \[\theta^i= \e^{-\frac{2b}{3}x}\ \hat\theta^i\] with 
\label{thetas}
\begin{eqnarray*}
\hat\theta^1&=&\d y -p\d x\\
\hat\theta^2&=&\d z-F\d x-2q(\d p-q\d x)\\
\hat\theta^3&=&-\tfrac{2^{4/3}}{\sqrt{3}}(\d p-q\d x)\\
\hat\theta^4&=&2^{-1/3}\d x\\
\hat\theta^5&=&
3 A_2(\d y-p\d x)+\tfrac{2^{2/3}}{3}b(\d p-q\d x)-2^{2/3}\d q +A_1
\d x,
\end{eqnarray*}
where 
\begin{eqnarray}
\label{a1}
A_1&=&
\tfrac{1}{2^{1/3}}
\left( a_1+2 a_2p+3 a_3 p^2+4 a_4 p^3+5 a_5p^4+6 a_6p^5+2bq\right),\\
\label{a2}
A_2&=&\tfrac{1}{45\cdot 2^{2/3}}\left( 9 a_2+27 a_3p+54 a_4 p^2+90  a_5p^3+ 135 a_6p^4+2b^2\right).
\end{eqnarray}
For further convenience we define 
\begin{eqnarray}
\label{a3}
A_3&=&
\tfrac{9}{20\cdot 2^{2/3}}
\left( a_3+4a_4p+10a_5p^2+20 a_6 p^3\right)
\\
\label{a4}
A_4&=&
\tfrac{9}{10}\left( a_4+5a_5p+15 a_6 p^2\right),
\label{a5}
 \\
 A_5
&=&
\tfrac{27}{4\cdot 2^{1/3}}
\left(a_5+6a_6p \right),\\
 A_6
&=&
\tfrac{243}{2\cdot  2^{2/3}}
a_6.
\label{a6}
\end{eqnarray}
Note that we use here a different metric in the conformal class $[g_F]$ than in \cite{nurowski07}. We have rescaled the metric in \cite{nurowski07} by $\e^{-\frac{4b}{3}x}$ which will give Cotton flat metrics for some $F$'s. When we write in the following ``not conformal'' we mean ``nowhere locally conformal''. Correspondingly, ``conformal'' for us always means ``locally conformal''.
\bs\label{a456not0}
If at least one of $a_4$, $a_5$, or $a_6$ is not zero, then the conformal class $[g_F]$ corresponding to $F=q^2+ \sum_{i=0}^6 a_i p^i+b z$ 
is not conformally Cotton and thus,  not conformally Einstein.
\es
\bprf

Recall that a metric which is  conformally Einstein is conformally Cotton. This  means that there exists a gradient field $T$ such that 
\begin{equation}\label{confcot}
C(T):=C+W(T,.,.,.)\equiv 0,
\end{equation} 
where $W$ is the Weyl tensor and $C$ is the Cotton tensor   (see e.g. \cite{gover/nurowski04}).
Writing $T=(\ups^1, \ldots, \ups^5)$  with $\ups^i=\theta^i(T)$ and
using the formulae in the appendix we get
$0\equiv C(T)_{112}= 
A_4\e^{\tfrac{4b}{3}x} \ups^4$.
Our assumption about $a_4$, $a_5$, and $a_6$ means that $A_4\not\equiv 0$. Thus, $\ups^4$ must be zero.
Furthermore, $0\equiv  C(T)_{214}= 
-A_4\e^{\tfrac{4b}{3}x} \ups^1$, which implies $\ups^1=0$. Finally, we get
$0\equiv C(T)_{314}= C_{314}=
-\tfrac{\sqrt{3}}{3} A_4\e^{2bx}\not=0$. 
 This means that with our assumptions  about $F$, the metric $g_F$ cannot be conformally Cotton, and hence, not conformally Einstein.\eprf

\bbem We observe the remarkable fact that  for any $F$ as in \eqref{eff}  the Riemann tensor of $\tg_F$   considered as an endomorphism of $\Lambda^2 T^*\tem$ has  rank   $\le 4$. In some cases  it can be even more degenerate. Hence, in order to obtain the $14$-dimensional group $\mathrm{G}_{2(2)}$ as holonomy group also derivatives of the curvature have to contribute to the holonomy algebra.
\ebem

\bs \label{a3not0}
For $F=q^2+a_3p^3+a_2p^2 +a_1p +a_0 +bz$ 
with $a_3\not=0$  the
 metric $g_F$ is the unique Cotton flat metric in $[g_F]$, but $g_F$ is not conformally Einstein. 
 \es
\bprf
The assumptions on $F$ imply that $A_4\equiv 0$ and $A_5\equiv 0$. By the formulae in the appendix this implies that the Cotton tensor of  $g_F$  is zero.
Now we find the most general vector $T$ such that $C(T)_{jkl} = W_{ijkl}\ups^i=0$.
The formulae for the Weyl tensor give that
\[\begin{array}{rcccl}
W_{i514}\ups^i &=&  W_{1514}\ups^1&=&\e^{\tfrac{4b}{3}x}A_3 \ups^1 
\\
W_{i115}\ups^i &=&  W_{4115}\ups^4&=&-\e^{\tfrac{4b}{3}x}A_3 \ups^4,
\end{array}
\]
which imply that $\ups^1=\ups^4=0$. Using $\ups^1=0$, we get
\[0=
W_{i414}\ups^i \ =\  W_{2414}\ups^2\ =\ \e^{\tfrac{4b}{3}x}A_3 \ups^2,\]
and thus $\ups^2=0$. Now the condition
\[0=
W_{i114}\ups^i \ =\  W_{3114}\ups^3+W_{5114}\ups^5\]
gives $\ups^5=\tfrac{2^{4/3}}{\sqrt{3}} b\ups^3$. This turns out to
solve all the remaining equations (\ref{confcot}). Hence, the most general $T$ solving \eqref{confcot} is given by
$\ups^i=f(0,0,1,0, \tfrac{2^{4/3}}{\sqrt{3}} b)$ with a smooth function $f$.
To define a scale $\ups $ such that $\e^{2\ups}g_F$ is
Einstein, this $T$ must be a gradient, which means that $\d \ups  =g(T,.)=f
\left(\tfrac{2^{4/3}}{\sqrt{3}}b\theta^1 +\theta^3\right)$. Thus in
such a case $\d (f\tau )=0$, where
$\tau=\left(\tfrac{2^{4/3}}{\sqrt{3}}b\theta^1 +\theta^3\right)$.
Calculating $\d\tau$ we get
\[0=\d(f\tau)\wedge \theta^1\wedge \theta^3 = f\d\tau \wedge \theta^1\wedge \theta^3
=
f \tfrac{2}{\sqrt{3}}\e^{\tfrac{2b}{3}x}
 \theta^1\wedge\ \theta^3\wedge \theta^4\wedge\theta^5.
 \]
But this implies that $f\equiv 0$. Hence, $g_F$ is the unique (up to a constant) Cotton flat  metric in $[g_F]$. The formulae for $\ro$ show that it is not Einstein. Thus, there is no Einstein metric in $[g_F]$. 
\eprf
We can summarise the results about whether $[g_F]$ contains an Einstein metric  in
\btheo\label{theorem0}
Let $F$ be given by $F=q^2+ \sum_{i=0}^6 a_i p^i+b z$ with at least one of $a_3$, $a_4$, $a_5$, $a_6$ not equal to zero. Then the conformal class $[g_F]$ 
does not contain an Einstein metric. If furthermore, $a_4=a_5=a_6=0$, then $g_F$ is  Cotton flat.
\etheo

Now we study the property of $[g_F]$ whether it contains a metric $g$ with the properties
\eqref{null-line} and \eqref{null-ric}, which were subject to Theorem \ref{theo-null}.

\blem
Let $(M,g)$ be pseudo-Riemannian manifold that admits a null line $L$. Then 
$Ric(X,.)=0$ for all $X\in L^\bot$ if and only if locally there is a vector field $K$  tangent to $ L$ and smooth function $\phi$ such that  $Ric = \phi g(K,.)\otimes g(K,.)$. 
Each of these properties implies that $(M,g)$ has vanishing scalar curvature and thus $\ro=\frac{1}{n-2}Ric$.
\elem
\bprf This is easily verified in a basis.\eprf
\blem
Let $(M,g)$ be pseudo-Riemannian manifold that admits a $\nabla^g$-parallel null line $L$ and satisfies the condition 
that $Ric(X,.)=0$ for all $X\in L^\bot$. Then
 the Weyl tensor $W$ of $g$ and hence of every metric in the conformal class of $g$ satisfies
\begin{eqnarray}
\label{one}
W(.,K,K,X)=0
\text{, for all $ K\in L$ and  $X\in L^\bot$}
\end{eqnarray}
\elem
\bprf
Since $L$ is parallel, the curvature $R$ of $g$ satisfies $R(U,V,K,X)=0$ for all $U,V\in TM$, $K\in L$ and $X\in L^\bot$.
Then the property $Ric (X,.)=0$ and hence $\ro (X,.)=0$ yields \eqref{one}.
\eprf

\btheo\label{theorem01}
Let $F$ be given by $F=q^2+ \sum_{i=0}^6 a_i p^i+b z$ with at least one of $a_3$, $a_4$, $a_5$, $a_6$ not equal to zero. Then the conformal class $[g_F]$ does not contain a metric $g$ with the properties \eqref{null-line} and \eqref{null-ric}.
\etheo

%


\bprf
We consider the most general null line $L$
 for $[g_F]$.
 We will show that there is no metric $g$ in the conformal class $[g_F]$ such that 
 conditions \eqref{null-line} and \eqref{null-ric} hold for $L$.
Let $K^i$ be tangent to $L$. We have the following four cases  to be excluded:
\bnum
\item
[a)]
 $K^i= (1, \alpha, \beta, \gamma, \alpha\gamma- \frac{1}{2} \beta^2)$,
\item
[b)] $K^i=( 0, 1, \beta , \frac{1}{2}\beta^2, \gamma)$,
\item
[c)] $K^i=( 0,0,0, 1, \gamma)$,
\item
[d)] $K^i=( 0,0,0, 0,1)$,
\enum
where $\alpha$, $\beta$, and $\gamma$ are arbitrary functions.
This is achieved by analysing the conformally invariant condition  \eqref{one} and, in cases b) and d), the properties 
\eqref{null-line} and \eqref{null-ric}, i.e.  $\hat{\nabla}_aK_b=f_aK_b$ and $\hat{P}_{ab}=\Phi K_a K_b$ for all metrics $\hat g\in [g_F]$.
The calculations, in which we will refer to the polynomials as defined in \eqref{a3} and \eqref{a4},
are based on the formulae provided in the appendix. Recall that  $A_3\equiv 0$ means that $g_F$ is conformally Einstein, and $A_4\equiv 0$ means that $g_F$ is conformally Cotton and not conformally Einstein if $a_3\not=0$. 
\begin{itemize}
\item
[c)]  Case c) is excluded because it is in contradiction with $g_F$ not being conformal to Einstein: One of the vectors from $K^\bot$ is $X^i= (1,\gamma, 0,0,0 )$. For this we get that $W_{2bcd}K^bK^d X^c=
A_3\e^{\tfrac{4b}{3}x } $. Hence,
 condition \eqref{one}  implies  that $g_F$ is conformal to an Einstein metric.
\item[a)]
First we exclude case a)  in the not conformally Cotton case, i.e. when at least one of $a_i\not=0$ for $i=4,5,6$, i.e. $A_4\not\equiv 0$.
In this case equation \eqref{one} for $X^i=(1,0,0,0 ,\einhalb\beta^2-\alpha \gamma)\in K^\bot$ gives
\[
0= W_{5bcd}K^bK^dX^c =
-A_3\gamma\e^{\tfrac{4b}{3}x},
\]
and thus $\gamma =0$. This yields 
\[
0=W_{1bcd}K^bK^dX^c = 
-2A_4\beta^2\e^{\tfrac{4b}{3}x},
\]
and therefore $\beta=0$.  Hence $X^i=( 1, 0,0,0,0)$. This gives 
\[0=
W_{4bcd}K^bK^dX^c = 
A_4\alpha\e^{\tfrac{4b}{3}x},
\]
and thus  $\alpha=0$, i.e. $K^i=(1,0,0,0,0)$.  Furthermore, for $Y^i=(0,1,0,0,0)\in K^\bot$ we get
\[0=W_{4bcd}K^bK^dY^c=
-A_4\e^{\tfrac{4b}{3}x}
\not=0,\] which gives the contradiction. 

Now we exclude case a) when $g_F$ is  Cotton flat, i.e. $a_4=a_5=a_6=0$, which means $A_4\equiv0$ and $a_3\not=0$. Equation \eqref{one} gives 
\[0=W_{5bcd}K^bK^dX^c= -
A_3\gamma \e^{\tfrac{4b}{3}x}\] which implies $\gamma=0$.
Furthermore, for  $Z^i=(0,0,1,0,-\beta)\in K^\bot$ we obtain from Eq. \eqref{one} that
\[
W_{4bcd}K^bK^dZ^c=
-\tfrac{A_3}{\sqrt{3}\cdot 2^{2/3}}
\e^{\tfrac{4b}{3}x}\left(
4b+
\sqrt{3}\cdot  2^{2/3}\beta
\right)
\]
This means that $\beta= -  \tfrac{2^{4/3}}{\sqrt{3}}b$.
Using this, equation \eqref{one} for $U^i=(0,0,0,1,\alpha)\in K^\bot$ gives
\[0=
W_{5bcd}K^bK^dU^c= A_3\e^{\tfrac{4b}{3}x},\] which is in contradiction with $A_3\not=0$.
\item[b)] 
For  case b) equation (\ref{one}) with vector $X^i=(1,0,0, \gamma,0)\in K^\bot$ gives 
\[W_{4bcd}K^bK^dX^c= -\tfrac{A_3}{2}\beta^2\e^{\tfrac{4b}{3}x}
=0,\] and thus $\beta=0$. 
Hence, in this case we have $K^i=(0,1,0,0,\gamma)$. 
Now we calculate $\hat\nabla K$ for the metric $\hat{g}=\mathrm{e}^{2\Upsilon}g$ with an arbitrary function $\Upsilon=\Upsilon(x, y, z , p, q)$. The condition that $K$ is tangent to a parallel null line for some $\Upsilon$ implies that the  first component of $\hat\nabla K$ must be zero. This implies 
that $\Upsilon=\Upsilon(x,y, p)$.
Using this we find that the third component of $\hat\nabla K$ vanishes if and only if   \[0\ =\ 
\tfrac{ \sqrt{3}}{2^{4/3}} \gamma \partial_p\Upsilon \theta^1 - \tfrac{\sqrt{3}}{3} (\gamma + 3\cdot 2^{-4/3} \partial_p\Upsilon )\theta^4 .\] This  yields $\Upsilon=\Upsilon(x,y)$ and $\gamma=0$.
With $\gamma=0$,  the fifth component of  $\hat\nabla K$ must vanish, which implies $\Upsilon=\Upsilon(x)$. Calculating the Schouten tensor $\ro$  for such $\Upsilon$ we find that $\hat \ro_{14}$ vanishes if and only if $A_3$ vanishes. Since in the conformally non-Einstein case the quantity $A_3$ is non-vanishing we get a contradiction with the condition $\hat \ro_{14}=0$ which is implied by $\hat \ro_{ab} = \phi K_a K_b$ and the $K$ with $\beta=\gamma=0$. This excludes the case b).
\item[d)]
A similar argument can be used in the case d). Here $K=(0,0,0,0,1)$ and the most general choice of the metric $\hat g=\mathrm{e}^{2\Upsilon}g$ leads to the following formula for the second component of the covariant derivative of $K$: $(\hat\nabla K)^2\wedge\theta^1\wedge\theta^2=-\tfrac{1}{2\sqrt{3}}\heta^3\wedge\heta^1\wedge\heta^2$. The condition that $K$ is parallel in $\hat\nabla$ requires that this must be zero, which excludes the case d) as well. 
\end{itemize}
Hence, for an arbitrary null vector and any metric in the conformal class $[g_F]$ we have shown that conditions \eqref{null-line} and \eqref{null-ric} cannot be satisfied together.\eprf

\section{Ambient metrics with holonomy $\mathrm{G}_{2(2)}$}\label{spinsec}
For those  conformal classes introduced in the previous  section that
are not conformally Einstein the relation between the holonomy of the
ambient metric and the conformal holonomy is more involved than in the
conformally  Einstein case.  We will now show that for some $[g_F]$ the ambient metric has holonomy exactly $\mathrm{G}_{2(2)}$.  The strategy is to show that the ambient manifold admits exactly one parallel spinor which is not null and exclude the existence of holonomy invariant null spaces by using Theorem \ref{theorem01}.

In Theorem \ref{maintheo} in the introduction we have proven that the ambient metric for $g_F$ as defined in \eqref{gF} is given as 
\[-2\d u\d t + t^2 g_F - 2tu \ro - u^2 B.\]
Calculating $\ro$ and $B$ explicitly in the appendix this reads as
\begin{eqnarray}\label{theoambient}
\tg_F&=&-~2 \d  t\d  u + t^2 g_F + 2tu~ \e^{\tfrac{4b}{3}x} \left( A_4(\theta^1)^2 + 2A_3 \theta^1\theta^4 +A_2 \theta^4 \right) +
\\
\nonumber
&&{ }+ \tfrac{1}{6}u^2~\e^{\tfrac{8b}{3}x} \left( A_6\theta^1\theta^2 + 2A_5 \theta^1\theta^4 +A_4 (\theta^4)^2 \right),
\end{eqnarray}
where the $A_i$'s are defined in (\ref{a2} -- \ref{a6}).
Note that the choice of a different $g_F$ in the conformal class than in \cite{nurowski07} results in a different coordinate system in which the ambient metric is expressed. Note also
 that in this form the ambient metric for $[g_F]$ has no $u^2$ terms  if $g_F$ is Cotton flat, i.e. if $a_4=a_5=a_6=0$. This means that for such $F$ it truncates at the  same order as the ambient metric of a conformal class with an Einstein metric, although it does not contain an Einstein metric if $a_3\not=0$.

In order to absorb the terms in the ambient metric coming from the terms of first and second order in $u$, we introduce the following co-frame on $M$:
\be
\eta^1&=&t \theta^1
\\
\eta^2&=&
 t \theta^2 
 - \tfrac{1}{12} \e^{\tfrac{4b}{3}x} \frac{u}{t}\left( 12A_2 t +A_4 \e^{\tfrac{4b}{3}x} \right)\theta^4
  \\
  &&{ }+
 \tfrac{1}{12} \e^{\tfrac{4b}{3}x} \frac{u}{t}
 \left( -24 A_3 t +12\cdot 2^{1/3} A_4pt -2 A_5 \e^{\tfrac{4b}{3}x}   u +2^{1/3} A_6 \e^{\tfrac{4b}{3}x} pu\right) \theta^1
\\
 \eta^3&=&
t \theta^3\\
 \eta^4&=&
 t\theta^4\\
 \eta^5&=&
 t \theta^5+
 +
 \tfrac{1}{12}\e^{\tfrac{4b}{3}x} \frac{u}{t}\left(
 12 A_4 +A_6 \e^{\tfrac{4b}{3}x}  u\right)\left( \theta^1 +2^{1/3}p\theta^4\right).
\ee
Then we write the ambient metric as 
\[ \wt{g}_F= - 2\d t\d u + 2\eta^1\eta^5 - 2\eta^2\eta^4 +(\eta^3)^2.
\]
For the calculation of the parallel spinor we use the orthonormal basis
\label{onframe}
\[
\begin{array}{llll}
\displaystyle \xi^0=\frac{1}{\sqrt{2}}( \d t-\d u),&
\displaystyle\xi^1=\frac{1}{\sqrt{2}}(  \eta^1+ \eta^5),&
\displaystyle\xi^2=\frac{1}{\sqrt{2}}(  \eta^2-  \eta^4),&
\displaystyle\xi^3=\eta^3\\[4mm] 
\displaystyle\xi^4=\frac{1}{\sqrt{2}}(  \eta^2+\eta^4),&
\displaystyle\xi^5=\frac{1}{\sqrt{2}}(  \eta^1- \eta^5),&
\displaystyle \xi^6=\frac{1}{\sqrt{2}}(\d t+\d u),&
\earr
in which $\tg_F$ reads as
\[\wt{g}_F\ =\ \tg_{ij} \xi^i\xi^j\ =\  (\xi^0)^2 +(\xi^1)^2+(\xi^2)^2 +(\xi^3)^2 -(\xi^4)^2 -(\xi^5)^2-(\xi^6)^2.
\]
We represent the Clifford algebra 
$\mathrm{Cl}(4,3)$ by means of $\s$-matrices satisfying the relation
\begin{equation}\label{cliff}
\s_i\s_j + \s_j\s_i = 2 \tg_{ij}\1_8.\end{equation}
They are given as:

\[\begin{array}{cccc}
\s_{0}= \left(\begin{array}{cc}
0&\gamma_0 \\ \gamma_0 & 0
\end{array}\right),&
\s_{1}= \left(\begin{array}{cc}
0&\gamma_2\\\gamma_2&0
\end{array}\right),
&
\s_2= \left(\begin{array}{cc}
0&\gamma_4\\\gamma_4&0
\end{array}\right)
&
\s_3= \left(\begin{array}{cc}
\1_4&0\\0&-\1_4
\end{array}\right)
\\[4mm]
\s_{4}= \left(\begin{array}{cc}
0&\gamma_1 \\ \gamma_1 & 0
\end{array}\right),&
\s_{5}= \left(\begin{array}{cc}
0&\gamma_3\\\gamma_3&0
\end{array}\right),
&
\s_6= \left(\begin{array}{cc}
0&-\1_4\\\1_4&0
\end{array}\right).
&
\end{array}\]
where
\[
\begin{array}{lll}
\gamma_{0}=
\left(\begin{array}{rrrr}
0&0&0&1\\
0&0&1&0\\
0&1&0&0\\
1&0&0&0
\end{array}\right),&
\gamma_2=
 \left(\begin{array}{rrrr}
0&0&1&0\\
0&0&0&-1\\
1&0&0&0\\
0&-1&0&0\\
\end{array}\right),&
\gamma_4=  \left(\begin{array}{rrrr}
1&0&0&0\\
0&1&0&0\\
0&0&-1&0\\
0&0&0&-1
\end{array}\right)
,\\[8mm]
\gamma_1=
 \left(\begin{array}{rrrr}
0&0&0&-1\\
0&0&1&0\\
0&-1&0&0\\
1&0&0&0
\end{array}\right),
&
\gamma_3= \left(\begin{array}{rrrr}
0&0&-1&0\\
0&0&0&-1\\
1&0&0&0\\
0&1&0&0\\
\end{array}\right).
\end{array}\]
Note that \[
\gamma_i^2=(-1)^i \1_4\]
which implies relation \eqref{cliff}. 
The invariant scalar product $\la.,.\ra$ is given by 
\[ 
\la \vf ,\psi\ra := -\left( \s_4 \cdot \s_5 \cdot\s_6 \cdot \vf,\psi\right),\]
where $(.,.)$ is the Euclidean standard scalar product on $\rr^8$. In the standard basis of $\rr^8$ the split signature scalar product $\la.,.\ra$ is given by the matrix
\[ 
\left(\begin{array}{rrrr}
0&0& \J_2&0\\
0&0&0&-\J_2\\
-\J_2&0&0&0\\
0&\J_2&0&0
\end{array}\right),
\]
where $\J_2 = 
\left(\begin{array}{rr}
0&-1\\ 1&0
\end{array}\right)$. It satisfies the relation 
\begin{equation}\label{realselfad}
\la \s_i\cdot \vf,\psi\ra = - \la \s_i\cdot\psi, \vf\ra,\end{equation}
 which implies its invariance.
Hence, the scalar product gives a metric on the spin bundle, which we denote by the same symbol, and which   is parallel w.r.t. the lift of the Levi-Civita connection $\widetilde{\nabla}$.
 
Then we have to solve the parallel spinor equations 
\begin{equation}\label{parspin}0=\widetilde{\nabla} \psi =\d \psi+\frac{1}{4} \sum_{k,l=0}^6 \widetilde{\Gamma}^{kl}\s_{k} \s_{l} \psi.\end{equation}
Here $\widetilde{\Gamma}^{kl}$ are the Levi-Civita connection $1$-forms for the ambient metric $\tg_F$ in the orthonormal co-frame $\xi^i$. I.e., $\wt{\Gamma}^{ij}$ are determined by
$\wt{\Gamma}^{ij}=-\wt{\Gamma}^{ji}$, $d\xi^i+\wt{\Gamma}^{i}_{\ j}\wedge \xi^j=0$, and $\wt{\Gamma}^{ij}=\wt{\Gamma}^{i}_{ \ k}g^{kj}$.

\bs\label{solution}
Let $F=q^2+ \sum_{i=0}^6 a_i p^i+b z$. 
Then the 
non-null spinor
\[\psi=
\left(
0,
-\e^{\tfrac{b}{3}x} ,
 \e^{\tfrac{b}{3}x} , 
0 , 
\sqrt{\tfrac{2}{3}} \e^{\tfrac{b}{3}x}\left(2^{\tfrac{1}{3}} b \e^{\tfrac{2b}{3}x}-3\right),
0,
0,
\sqrt{\tfrac{2}{3}} \e^{\tfrac{b}{3}x}\left(2^{\tfrac{1}{3}} b \e^{\tfrac{2b}{3}x}+3\right)
\right)
\]
 is a solution of the parallel spinor equation \eqref{parspin}. In particular, the holonomy of the ambient metric of $[g_F]$ is contained in $\mathrm{G}_{2(2)}$.
 \es
 \begin{proof}
 One checks by direct calculations that $\psi$ is parallel and not null with $\la \psi,\psi\ra=4\sqrt{6}$.
\eprf
 
For  completeness we will give below a formula for the parallel three-form $\w$ that defines the $\mathrm{G}_{2(2)}$ structure. The form $\w$ is related to the spinor $\psi$ by the following relation (see for example \cite{kath98g2}): First one defines a skew $(2,1)$-tensor $A^\psi$ depending on $\psi$ via
\[ X\cdot Y \cdot \psi - \wt{g}(X,Y)\psi =A^\psi(X,Y)\cdot \psi\]
and obtains $\w$ by dualising it
\[\w (X,Y,Z):= \wt{g} (X,A^\psi(Y,Z)).\]
Calulating this with Mathematica we get that $\w$ is equal to
\be
\w&=&
\tfrac{1}{6\cdot 2^{5/6} \sqrt{3}}\left( 18 f(-x)-3\cdot 2^{1/3} f(x) +4b^2 f(x)\right)\left( \xi^{012}-\xi^{146}\right)
\\
&&
+\tfrac{1}{6\cdot 2^{5/6} \sqrt{3}}\left( 18 f(-x)+3\cdot 2^{1/3} f(x) -4b^2 f(x)\right) \left(\xi^{014}+\xi^{126}\right)
\\
&&
-\tfrac{1}{6\cdot 2^{5/6} \sqrt{3}}\left( 18 f(-x)+3\cdot 2^{1/3} f(x) +4b^2 f(x)\right) \left(\xi^{025}
-\xi^{456}\right)
\\
&&
+\tfrac{1}{6\cdot 2^{5/6}\sqrt{3}}\left( -18 f(-x)+3\cdot 2^{1/3} f(x) +4b^2 f(x)\right)\left( \xi^{045}+\xi^{256}\right)
\\
&&
-\tfrac{ 2^{1/6} b \left( 3\cdot 2^{2/3} - bf(x)\right)}{\sqrt{3}\left(-3+2^{1/3} b f(x)\right) }\left(\xi^{016}
+\xi^{124}\right)
+ \tfrac{2^{1/3b}}{3}f(x) \left(-\xi^{023}-\xi^{034}+\xi^{236} +\xi^{346}\right)
\\
&&
 +\xi^{036} +\xi^{135}+\xi^{234} +\tfrac{2^{5/6}b}{\sqrt{3}}\left(\xi^{056}-\xi^{245}\right)
\ee
with
\[
f(x):=\e^{\tfrac{2b}{3}x}\]
and $\xi^{ijk}:=\xi^i\wedge\xi^j\wedge\xi^k$, where $\xi^i$ is the orthonormal coframe given on page \pageref{onframe}. 
A direct calculation verifies that $\w$ and its Hodge dual are closed.

In order to conclude the proof of Theorem \ref{maintheo} it only remains to prove the following lemma.

\blem \label{spinlemma} If a $7$-dimensional spin manifold $\tem$ with metric $\tg$ of signature $(4,3)$ admits a parallel non null spinor $\psi$  and a parallel line of null spinors, then there is a parallel line bundle of tangent vectors on $\tem$.
\elem
\bprf
We fix a spinor $\vf$ that spans the parallel line of spinors. There is a 1-form $f$ such that $\tnab
\vf=f\otimes \vf$.  We associate to $\psi$ and $\vf$ a vector field $V$ via transposing the Clifford multiplication, i.e.
 \[\tg(V,X) \ =\  \la X\cdot \psi, \vf\ra \ =\ -\la X\cdot \vf,\psi\ra
 \]
 for all $X\in T\tem$.
 The well known formula 
 \[
Y\left(\la X\cdot \vf, \psi\ra\right) =
\la \tnab_Y X\cdot \psi,\vf\ra + \la X\cdot \tnab_Y\psi,\vf\ra +\la X\cdot \psi, \tnab_Y\vf\ra
\]
for two spinor fields $\vf$ and $\psi$, and two vector fields $X$ and $Y$, shows that 
 $V$ spans a parallel line. Indeed,
 it implies
 that
 \be
 \tg(\tnab_XV,Y)&=&X(\tg(V,Y))-\tg(V,\tnab_XY)\\
 &=& X(\la Y\cdot \psi,\vf\ra)-\la \tnab_X Y\cdot \psi,\vf\ra\\
 &=& f(X) \la Y\cdot \psi, \vf \ra \\
 & =& f(X) \tg (V,Y),\ee
 for all $X,Y\in T\tem$. 
For the proof, 
we have to exclude that $V\equiv 0$ i.e. that
\begin{equation}
\label{null}g(V,X) \ =\  \la X\cdot \psi, \vf\ra =0\end{equation}
for all $X\in T\tem$.
We will show that
this contradicts   $\psi$ being not null and  $\vf$ being null.
To this end,  
at each tangent space $T_p\tem =\rr^{4,3}$, 
 consider the map 
 \[
\rr^{4,3}\ni X\mapsto X\cdot \psi\in \D_{4,3}.\]
 Using the transitive action of $\mathrm{Spin}(4,3)$ on spheres in $\rr^{4,4}$, one  shows \cite{kath98g2} that $\psi$ is not null if and only if this map has a trivial kernel. Hence, with $\psi$ being not null, the vector space 
 \[W:=\{X\cdot \psi\mid X\in \rr^{4,3}\}\]
 has dimension seven.
Furthermore, property \eqref{realselfad} implies  
\begin{equation}\label{vau}
 2 \la X\cdot \psi, Y\cdot \psi\ra = -g(X,Y)\la \psi,\psi\ra,
\end{equation}
for all $X,Y\in \rr^{4,3}$.
Since $\psi$ is not null, this shows
that 
$W\subset \rr^{4,4}$ is non-degenerate. Equation \eqref{null} then implies that $\rr \vf= W^\bot$ which is contradicts $\vf$ being null.
\eprf
This lemma shows that the existence of a parallel maximal totally null subspace  yields the existence of a parallel line bundle in the tangent bundle for the ambient metric, which by Theorem \ref{einsteintheo} contradicts $g_F$ not being conformally Einstein.

\bbem
Note also that  the existence of  two non null spinor yields the existence of a parallel vector field. This is true by the result in \cite{kath98g2} that the isotropy group of two spinors that are not null is given by $\mathrm{SU}(1,2)$ or $\mathrm{SL}(3,\rr)$. Both cases imply that there is a parallel vector field on $M$ that is not null.  
\ebem

\section*{Appendix}
Here we will give formulae for the Levi-Civita connection of $g_F$, its Schouten, Weyl, Cotton and Bach tensor.
$g_F$ is given as in \eqref{gF}, with $\theta^i$'s as on page \pageref{thetas}.
In this coframe the Levi-Civita connection 1-forms, i.e. matrix-valued
1-forms satisfying 
$\d \theta^\mu+\Gamma^\mu_{~\nu}\wedge\theta^\nu=0$,
$\Gamma_{\mu\nu}+\Gamma_{\nu\mu}=0$,
$\Gamma_{\mu\nu}=g_{\mu\sigma}\Gamma^\sigma_{~\nu}$, are: 
\begin{eqnarray*}
\Gamma_{12}&=& \Gamma_{23}\ =\ \Gamma_{25}\ =\ 0
\\
\Gamma_{34}&=&  - \tfrac{2^{1/3}}{3}b\hat\theta^3 +\tfrac{1}{\sqrt{3}} \hat\theta^5
\\
\Gamma_{35}&=& -\tfrac{1}{\sqrt{3}} \hat\theta^4
\\
\Gamma_{45}&=& \tfrac{2^{1/3}}{3}b\hat\theta^1 +\tfrac{1}{2\sqrt{3}} \hat\theta^3
\\
\Gamma_{15}&=& -\tfrac{2^{1/3}}{3}b \hat\theta^4
\\
\Gamma_{24}&=& -\tfrac{2^{1/3}}{3} b\hat\theta^4
\\
\Gamma_{13}&=&  -2\sqrt{3}\left( A_3\hat\theta^1+ A_2\hat\theta^4\right)
\\
\Gamma_{14} &=& 
2^{4/3}\left( 2^{4/3}A_3 q - A_2 b\right) \heta^1
+\tfrac{3\sqrt{3}}{2}A_2\hat\theta^3
-
\tfrac{2^{1/3}}{3}b\hat\theta^5
\end{eqnarray*}
where $A_1$ and $A_2 $ are defined in \eqref{a1}, \eqref{a2} and \eqref{a3}.
 Then the Schouten tensor is given as
 \[
 \ro= 
 -A_4(\heta^1)^2 -2 A_3\heta^1\heta^4 -A_2(\heta^4)^2
 \]
 with $A_4$ defined in \eqref{a4}. 
Let $W_{ijkl}$ be the Weyl tensor and  $W_{ij}$ be the 2-forms defined by $W_{ij}=\einhalb W_{ijkl}\theta^k\theta^l$. They are are given by
\begin{eqnarray*}
W_{12}&=& -A_4 \heta^1\wedge \heta^4
\\
W_{13}&=& -2A_4 \heta^1\wedge \heta^3 + \tfrac{2^{4/3}}{\sqrt{3}}
\left( 3\cdot 2^{1/3} A_4 q-A_3 b\right)\heta^1\wedge\heta^4
\\
W_{14}&=& -A_4 \heta^1\wedge \heta^2 + \tfrac{2^{4/3}}{\sqrt{3}}\left( 3\cdot 2^{1/3} A_4 q-A_3b\right)\heta^1\wedge\heta^3
+ \\
&&{ }+
\tfrac{1}{3}\left( 27A_2^2 - 12 \cdot 2^{1/3} A_1A_3 -6\cdot 2^{2/3}A_2 b^2 +40 A_3bq -24\cdot 2^{1/3} A_4q^2\right)
\heta^1\wedge \heta^4\\
&&{ }
+A_3 \left( \heta^1\wedge\heta^5 +\heta^2\wedge\heta^4\right)
\\
W_{15}&=&W_{24}\ =\ A_3\heta^1\wedge\heta^4
\\
W_{23}& =& W_{25}\ =\ W_{34}\ =\ W_{35}\ =\ W_{45}\ =\ 0.
\end{eqnarray*}
If $C_{ijk}$ is the Cotton tensor, the 1-forms $C_i=\einhalb C_{ijk}\theta^j\wedge\theta^k$ are given by
\begin{eqnarray*}
C_2&=& C_5\ =\ 0\\
C_3&=& 
-\tfrac{\sqrt{3}}{3}  A_4 \e^{2bx}\theta^1\wedge\theta^4
\\
C_4
&=& 
A_4 \e^{2bx}\left(  -\tfrac{\sqrt{3}}{3} \theta^1\wedge\theta^3
+
\tfrac{2^{2/3}}{3} q \theta^1\wedge\theta^4\right)
\\
C_1&=&
-\tfrac{\sqrt{3}}{3}A_5 \e^{2bx}
 \theta^1\wedge \theta^3 
 +
 \tfrac{2^{1/3}}{3}
\left( A_4 b + 2^{4/3} A_5 q \right)\e^{2bx}\theta^1\wedge\theta^4,
\end{eqnarray*}  
where $A_5$ is defined in \eqref{a5}.
Finally, the Bach tensor is given by
\[
B\ =\ -\tfrac{1}{6}~\e^{\tfrac{8b}{3}x} \left( A_6\theta^1\theta^2 + 2A_5 \theta^1\theta^4 +A_4 (\theta^4)^2 \right).
\]
These formulae enable the reader to calculate the connection coefficients of the truncated ambient metric
\[\tg_F=-2\d u\d t + t^2 g_F - 2tu \ro - u^2 B.\]
\def\cprime{$'$}

\end{document}